# Topology optimization on complex surfaces based on the moving morphable component (MMC) method and computational conformal mapping (CCM)

Wendong Huo[1], Chang Liu[1,2]*, Zongliang Du[1,2], Xudong Jiang[1], Zhenyu Liu[3], Xu Guo[1,2]*

*[1]State Key Laboratory of Structural Analysis for Industrial Equipment,*
*International Research Center for Computational Mechanics,*
*Dalian University of Technology, Dalian, 116023, P.R. China*

*[2]Ningbo Institute of Dalian University of Technology, Ningbo, 315016, P.R. China*

*[3]Changchun Institute of Optics, Fine Mechanics and Physics (CIOMP),*
*Chinese Academy of Sciences, Changchun 130033, P.R. China*

**Abstract**

In the present paper, an integrated paradigm for topology optimization on complex surfaces with arbitrary genus is proposed. The approach is constructed based on the two-dimensional (2D) Moving Morphable Component (MMC) framework, where a set of structural components are used as the basic units of optimization, and computational conformal mapping (CCM) technique, with which a complex surface represented by an unstructured triangular mesh can be mapped into a set of regular 2D parameter domains numerically. A multi-patch stitching scheme is also developed to achieve an MMC-friendly global parameterization through a number of local parameterizations. Numerical examples including a saddle-shaped shell, a torus-shape shell and a tee-branch pipe are solved to demonstrate the validity and efficiency of the proposed approach. It is found that compared with traditional approaches for topology optimization on 2D surfaces, optimized designs with clear load transmission paths can be obtained with much fewer numbers of design variables and degrees of freedom for finite element analysis (FEA) via the proposed approach.

*Corresponding authors. E-mail: c.liu@dlut.edu.cn (Chang Liu), guoxu@dlut.edu.cn (Xu Guo)

# 1. Introduction

Topology optimization aims at distributing a certain amount of material in a prescribed design domain in order to satisfy some design requirements and at the same time achieve exceptional performances. As a revolutionary and powerful design method, it can help engineers create competitive designs in a systematic way and has attracted the attention of many engineering fields since its invention. After the pioneering work of Bendsøe and Kikuchi [1], numerous topology optimization approaches such as the Solid Isotropic Material with Penalization (SIMP) method (also named as variable density method) [2, 3], evolutionary structural optimization (ESO) method [4], level set method [5, 6], just to name a few, have been established and applied successfully in many engineering applications.

Shell, as a typical engineering structure surrounding space in an aesthetic way, enjoys the benefits of efficient load-carrying capacity and high stiffness [7]. Topology optimization of shell structures could further promote their strength to weight ratio and help achieve a lightweight structural design in many application fields such as mechanical, civil, marine, and aeronautical engineering. From the mathematical point of view, topology optimization of shell structures is equivalent to finding the optimized material distribution on a 2D manifold (roughly speaking, a region that can be parameterized by two parameters). Compared to topology optimization on flat 2D and 3D space, corresponding research works on 2D manifold-based topology optimization are relatively rare. Most of the approaches are constructed under the implicit variable density framework. For example, Gea et al. [8] investigated the optimal design of 3D plate and shell structures for both static and dynamic cases using the density method. The filter scheme of the SIMP method on surfaces and topology optimization on two-dimensional manifolds was proposed for several physical fields by Deng et al. [9]. Sigmund and his co-authors combined high-performance computing and topology optimization to design ultra-large-scale shell structures [10]. Other recent excellent progress on topology optimization of shell structures can be found in [11–14] and the references therein.

The challenging issues associated with the SIMP-based approach mainly come from two sources. One is that compared with its flat 2D counterpart, topology optimization on curved surfaces involves larger numbers of design variables and degrees of freedom for optimization and finite element analysis (FEA), respectively. This is due to the fact that more finite elements are required to discretize a shell-type structure in order to guarantee the accuracy of both geometry modeling and FEA, especially when the shell is of complex shape and has relatively large local curvatures. The other is that the filter approaches, which are very effective in suppressing the numerical instabilities (e.g., checkerboard pattern) in flat 2D case, may encounter some difficulties when applied to unstructured meshes generated on curved surfaces since it is difficult to determine the (element-wise) radius of the filter in advance, which should be local curvature-dependent under the considered case.

Topology optimization of shell structures has also been investigated with the level set approach. One of the excellent works is presented by Ye et al. [15], in which the conformal geometry theory was first introduced to the field of structural topology optimization. By constructing the parametric domain of surfaces through conformal mappings, the classical level set topology optimization method is extended to manifolds. Since global parameterization of the manifold is pursued, the relatively high nonlinear mapping would increase the difficulties in the corresponding solution process.

It should be noted that the above works of topology optimization on surfaces mainly use implicit methods and therefore may suffer from problems such as grey elements, islanding effect, and an enormous number of design variables. Recently, explicit approaches have received more and more attention in the field of topology optimization [16-18]. Among them, the Moving Morphable Component (MMC) method describes the optimized structure using a set of geometrically explicit components. By taking its advantages of explicit description and a fewer number of design variables, the MMC method has been successfully extended to consider manufacturability [19], geometrical nonlinearity [20], dynamic performance [21], and multi-physics effects [22, 23], etc. Nevertheless, the current MMC method is developed on 2D or 3D Euclidean space and cannot be directly applied to design optimization on manifolds effectively.

Compared to topology optimization in two-dimensional flat space, it is necessary to solve the following challenging problems when the MMC approach is developed for achieving topological design on a complex surface $\mathcal{S}$ (in general can be considered as a two-dimensional manifold) with an arbitrary genus. Firstly, how to describe the geometry configuration of a complex surface in a universal and flexible way? Secondly, how to construct the topological description function (TDF) $\phi$ of a component whose support set is lying entirely on the surface (i.e., $\mathrm{Supp}\ \phi \subset \mathcal{S}$)? Thirdly, how to carry out MMC-based topology optimization on a complex surface with a high genus number?

The present work intends to solve the aforementioned problems in an integrated way. An unstructured triangular mesh which is highly robust and flexible for geometry/topology description is employed to describe the embedding information of a smooth surface in $\mathbb{R}^3$ with arbitrary accuracy. With the use of computational conformal mapping (CCM) on triangular meshes, the traditional MMC approach originally established in flat space is extended to a simply-connected open surface with genus zero (which is homeomorphic to a planar rectangle unit cell) at first, and then generalized to account for arbitrary complex surfaces with the help of the multi-patch stitching scheme, with which an MMC-friendly global parameterization of a complex surface can be achieved through a set of local parameterizations.

The remainder of the article is organized as follows. The theoretical foundations of the proposed approach including the 2D MMC-based framework for topology optimization and computational conformal mapping (CCM) used for surface parameterization are introduced in Section 2. In Section 3, firstly, the mathematical formulation of the considered problem is provided. Secondly, taking a simply-connected open surface with zero genus as an example, the numerical algorithm for carrying out explicit topology optimization on a 2D manifold is described. Finally, the flowchart of the proposed integrated solution procedure applicable to complex surfaces with arbitrary genus is described in detail in the last part of this section. Three numerical examples are then investigated in Section 4 to demonstrate the effectiveness of the proposed approach. Some concluding remarks including the summary of the present work and discussions on possible directions of future researches are provided in the last section.

## 2. Theoretical foundation

### 2.1 Moving Morphable Component (MMC) method

The MMC method for topology optimization was developed in [16] to optimize structural topology in an explicit way. In this method, the basic units of optimization are a set of structural components and the variation of structural topology can be achieved by the moving, deforming, overlapping, and merging of these components. A typical 2D component (also adopted in the present work) is shown in Fig. 1 and its topology description function (TDF) can be described as

$$\phi = 1 - \left(\left(\frac{x'}{L}\right)^6 + \left(\frac{y'}{f(x')}\right)^6\right)^{1/6} \tag{2.1}$$

with

$$\begin{pmatrix} x' \\ y' \end{pmatrix} = \begin{bmatrix} \cos\theta & \sin\theta \\ -\sin\theta & \cos\theta \end{bmatrix} \begin{pmatrix} x - x_0 \\ y - y_0 \end{pmatrix} \tag{2.2}$$

$$f(x') = \frac{t^1 + t^2 - 2t^3}{2L^2}(x')^2 + \frac{t^2 - t^1}{2L}x' + t^3 \tag{2.3}$$

where $x_0$, $y_0$ and $\theta$ denote the two coordinates of the component's central point and its rotational angle with respect to the global Cartesian coordinate system, as illustrated in Fig. 1, respectively. In Eq. (2.1)-Eq. (2.3), the symbols $L$, $t^1$, $t^2$ and $t^3$ denote the length and three thickness parameters of the component, respectively.

Assuming that there are totally $n$ components existing in the design domain, the topology of the structure can be described by its TDF $\phi^s$ as

$$\begin{cases} \phi^s(x) > 0, & \text{if } x \in \Omega_s \\ \phi^s(x) = 0, & \text{if } x \in \partial\Omega_s \\ \phi^s(x) < 0, & \text{if } x \in D\backslash(\Omega_s \cup \partial\Omega_s) \end{cases} \tag{2.4}$$

where D represents the design domain and $\Omega_s$ is the region occupied by the structure. The TDF $\phi^s$ can be constructed by the TDF of each component in terms of K-S function as [24]:

$$\phi^s = \mathcal{KS}(\phi^1, \phi^2, \dots, \phi^n) = \left(\ln\left(\sum_{i=1}^{n}\exp(\zeta\phi^i)\right)\right)/\zeta \tag{2.5}$$

where $\zeta$ is a large positive number, e.g., $\zeta = 100$ and $\phi^i$ is the TDF of the $i$-th component. Under this circumstance, the vector of the design variables associated with the $i$-th component can be identified as $\boldsymbol{D}_i = \left(x_0^i, y_0^i, \theta_i, L_i, t_i^1, t_i^2, t_i^3\right)^\mathsf{T}$ and the vector of the design variables of the whole structure is $\boldsymbol{D} = (\boldsymbol{D}_1^\mathsf{T}, \boldsymbol{D}_2^\mathsf{T}, \dots, \boldsymbol{D}_n^\mathsf{T})^\mathsf{T}$.

## 2.2 Computational conformal mapping (CCM)

Roughly speaking, a surface is a 2D manifold with intrinsic non-zero curvature embedding in the three-dimension (flat) Euclidean space. The original MMC approach described in subsection 2.1 is developed in the 2D (or 3D) flat space and therefore cannot be applied directly to solve topology optimization problems on surfaces. The key problem is how to construct the global TDF $\phi_S^s(x)$ of the structure on the 2D surface manifold. A natural idea is to establish a homeomorphic mapping $f: S \to M$ between the concerned surface $S$ embedded in a flat 3D space and a 2D planar parameter domain $M$ (parameterization), construct the corresponding topology description function $\phi_M^s(p)$ on $M$ and then use the inverse mapping $f^{-1}: M \to S$ to obtain the topology description function $\phi_M^s(p)$ through $\phi_S^s(x) = \phi_M^s(p = f(x))$ (see Fig. 2 for reference).

Parameterizing a complex surface with an arbitrary genus is, however, not a trivial task. Fortunately, thanks to the development of computational geometry, powerful tools such as conformal/quasi-conformal mapping techniques have been established and applied successfully in many interesting applications [25, 26]. Specifically, according to the theorems provided in [27, 28], a complex surface with an arbitrary genus can be mapped to a simply-connected planar parameter domain through some appropriate homeomorphic mapping, which can be determined numerically by solving a series of partial differential equations (PDEs) [29–32]. In the present work, the computational conformal mapping algorithm developed in [32–35] is adopted to construct the homeomorphism used for parameterization. The corresponding solution procedure is described briefly as follows.

The present conformal mapping is composed of two quasi-conformal mappings. For the sake of simplicity, taking a simply connected oriented open surface $S$ (a two-dimensional manifold) with genus zero as an example, the first quasi-conformal mapping $h: S \to \mathbb{D} \subset \mathbb{C}$ establishing a topological homeomorphism between $S$ and a planar unit disk $\mathbb{D}$ in complex plane $\mathbb{C}$ is constructed by solving the following partial differential equation:

$$\begin{cases} \Delta_S h = 0, & \text{on } S, \\ h(\partial S) = \partial \mathbb{D} \end{cases} \quad (2.6)$$

using finite element method [36]. In Eq. (2.6), $\Delta_S$ represents the Laplace-Beltrami operator defined on surface $S$. Once $h = h(S)$ is determined, the second complex quasi-conformal mapping $g(z = x + iy) = u(x, y) + iv(x, y): \mathbb{D} \subset \mathbb{C} \to M \subset \mathbb{C}$ from $\mathbb{D}$ to a planar parametric domain $M$ (a standard rectangle in the present work) is constructed by finding the solutions of the following generalized Laplace equations [33, 34]:

$$\begin{cases} \nabla \cdot (\mathbf{A}(\nabla u)) = 0, & (2.7a) \\ \nabla \cdot (\mathbf{A}(\nabla v)) = 0, & (2.7b) \end{cases}$$

where $u = u(x,y)$, $v = v(x,y)$ and $\nabla(\cdot) = \left(\frac{\partial(\cdot)}{\partial x}, \frac{\partial(\cdot)}{\partial y}\right)^T$. Furthermore, $\mathbf{A} = \begin{pmatrix} \alpha_1 & \alpha_2 \\ \alpha_2 & \alpha_3 \end{pmatrix}$ and $\alpha_1 = \frac{(\rho-1)^2+\tau^2}{1-\rho^2-\tau^2}$, $\alpha_2 = -\frac{2\tau}{1-\rho^2-\tau^2}$, $\alpha_3 = \frac{(\rho+1)^2+\tau^2}{1-\rho^2-\tau^2}$ with $\rho + i\tau = \mu_{h^{-1}}$ denoting the Beltrami efficient of the mapping of $h^{-1}: \mathbb{D} \to \mathcal{S}$ (the inverse mapping of $h$) [33, 34].

On condition that the two quasi-conformal mappings $h$ and $g$ are determined, the required conformal mapping $f: \mathcal{S} \to \mathcal{M}$ can be constructed as $f = g \circ h$. For the case where the concerned surface $\mathcal{S}$ has complex topology and non-zero genus, we refer the readers to [32–35] for more technical details on constructing the conformal mapping function.

**2.3 TDF construction on complex surfaces**

*2.3.1 Simply-connected open surface with zero genus*

This is the simplest case for TDF construction. Once the parameterization from $\mathcal{S}$ to $\mathcal{M}$ expressed in terms of $f$ is constructed, we can first define the TDF on $\mathcal{M}$ obtaining $\phi_{\mathcal{M}}^S(p), p \in \mathcal{M}$ and then determine the TDF on $\mathcal{S}$ through $\phi_{\mathcal{S}}^S(x) = \phi_{\mathcal{M}}^S(p)$ with $p = f(x \in \mathcal{S})$ (see Fig. 3 for reference). This treatment is well-posed since $f$ represents a topological homeomorphism between $\mathcal{S}$ and $\mathcal{M}$.

*2.3.2 Complex surface with arbitrary genus*

When the topology of the concerned surface $\mathcal{S}$ is a complex manifold with a high genus, there is no homeomorphism between $\mathcal{S}$ and a planar rectangle $\mathcal{M}$. Under this circumstance, cutting operation should be used to generate a simply-connected open intermediate surface $\mathcal{S}^*$, which can be made topologically equivalent to a planar rectangle $\mathcal{M}$ by the conformal mapping technique described above. Fig. 4 (a) demonstrates the procedure of the cutting operation for a torus surface. Therefore, we can define the TDF on $\mathcal{S}^*$ as $\phi_{\mathcal{S}^*}^S(x) = \phi_{\mathcal{M}}^S(f^*(x))$ where $f^*: \mathcal{S}^* \to \mathcal{M}$ (see Fig. 4 (b) for reference). Considering the fact that in general $\phi_{\mathcal{S}^*}^S(x)$ may take different values on different sides of a specific cutting line $\Gamma_i$ on $\mathcal{S}$ ($\Gamma_i'$ and $\Gamma_i''$ on $\mathcal{S}^*$), we propose to define the value of $\phi_{\mathcal{S}}^S(x)$ in terms of $\phi_{\mathcal{S}}^S(x)$ as

$$\phi^{\mathrm{S}}_{\mathcal{S}}(x) = \begin{cases} \phi^{\mathrm{S}}_{\mathcal{M}}(f^*(x)), & \text{if } x \in \mathcal{S}\setminus\Gamma_i, \\ \mathcal{KS}\left(\phi^{\mathrm{S}}_{\mathcal{M}}(f^*(x')), \phi^{\mathrm{S}}_{\mathcal{M}}(f^*(x''))\right), & \text{if } x \in \Gamma_i, x' \in \Gamma'_i \text{ and } x'' \in \Gamma''_i \end{cases}$$
(2.8a)
(2.8b)

where $f^*(x')$ and $f^*(x'')$ denote the values of $f^*$ taking on $x' \in \Gamma'_i$ and $x'' \in \Gamma''_i$ respectively.

*2.3.3   Constructing TDF with multi-patch stitching approach*

Although, the global parameterization technique can be applied for any surface that has manifold property in principle, the constructed global mapping function $f$ may be too stiff, and it has the potential to induce some numerical instabilities when $\phi^{\mathrm{S}}_{\mathcal{S}}(x)$ is generated based on $f$. This issue can be resolved by the so-called multi-patch stitching approach suggested in the literatures [37–40]. This approach can reduce the distortion of parametrization and thus greatly enhance the fidelity of geometry description of components especially when the concerned surface has large local curvature, high genus and/or non-manifold properties.

In this approach, as shown in Fig. 5 (a), the surface $\mathcal{S}$ is decomposed onto $N_{\mathcal{U}}$ parts, i.e., $\mathcal{S} = \bigcup_{k=1}^{N_{\mathcal{U}}} \mathcal{U}_k$ and in general $\mathcal{U}_k \cap \mathcal{U}_l \neq \emptyset$, $k, l = 1, \dots, N_{\mathcal{U}}$. For each part $\mathcal{U}_k$, we can establish a conformal parameterization through a mapping $f_k: \mathcal{U}_k \to \mathcal{M}_k$ where $\mathcal{M}_k$ is a rectangle on the parametric domain. For the purpose of constructing the global TDF for $\mathcal{S}$, we can place a set of components on each $\mathcal{M}_k$ and then generate the TDF on $\mathcal{U}_k$ through $\phi^{\mathrm{S}}_{\mathcal{U}_k}(x) = \phi^{\mathrm{S}}_{\mathcal{M}_k}(f_k(x))$ with $x \in \mathcal{U}_k$. Considering the fact the intersection between two parts may be non-empty (i.e., $\mathcal{U}_k \cap \mathcal{U}_l \neq \emptyset$) which is sometimes necessary for rendering smooth connections of the components located on neighboring parts, the global TDF for $\mathcal{S}$, i.e., $\phi^{\mathrm{S}}_{\mathcal{S}} = \phi^{\mathrm{S}}_{\mathcal{S}}(x)$ is determined as (see Fig. 5 (b) for reference):

$$\phi^{\mathrm{S}}_{\mathcal{S}}(x|x \in \mathcal{U}_i) = \begin{cases} \phi^{\mathrm{S}}_{\mathcal{U}_i}(x), & \text{if } x \in \mathcal{U}_i \setminus \bigcup_{k=1, k\neq i}^{N_{\mathcal{U}}} \mathcal{U}_k, \\ \mathcal{KS}\left(\phi^{\mathrm{S}}_{\mathcal{U}_i}(x), \dots, \phi^{\mathrm{S}}_{\mathcal{U}_j}(x)\right), & \text{if } x \in \mathcal{U}_i \cap \dots \cap \mathcal{U}_j. \end{cases}$$
(2.9a)
(2.9b)

Numerical examples presented in Section 4 show that this approach is very effective on stabilizing the optimization process and guaranteeing the smooth transition of the components in the optimized designs. It is also worth noting that $\mathcal{U}_k$ may also be multi-connected and has a high genus. Under this circumstance, the cutting operation described in the last subsection is also applicable for constructing $\phi^{\mathrm{S}}_{\mathcal{U}_k}(x)$.

## 3. The statement of the problem and its solution procedure

This section is devoted to the description of the problem statement and the solution procedure of the considered problem.

### 3.1 Problem statement

In the present work, as shown in Fig. 6, compliance minimization by distributing a certain amount of isotropic linear elastic material (the upper bound of the volume fraction is $\overline{V}$) on a 3D region $\mathcal{B}$ with a small uniform thickness (i.e., $t \ll a$, $t \ll b$ with $t$, $a$ and $b$ denoting the thickness and the characteristic length scales of the other two directions of $\mathcal{B}$, respectively) is considered. Under this circumstance, shell model can be used for carrying out structural response analysis in a more efficient way (compared to the treatment where 3D elasticity theory is adopted). It is also assumed that $\mathcal{B}$ can be parameterized by a bijective mapping $\overline{\varphi}$ from a parametric domain $\Omega = \omega \times \left(-\frac{t}{2}, \frac{t}{2}\right) = \left\{\xi = (\xi^1, \xi^2, \xi^3) \middle| (\xi^1, \xi^2) \in \omega, \xi^3 \in \left(-\frac{t}{2}, \frac{t}{2}\right)\right\}$ such that $\mathcal{B} = \overline{\varphi}(\Omega)$ and $\mathcal{S} = \overline{\varphi}(\omega \times (\xi^3 = 0))$ is the mid-surface of $\mathcal{B}$. Therefore, the corresponding topology optimization problem formulation can be written as:

$$\text{Find } \boldsymbol{U} = \boldsymbol{U}(\boldsymbol{u}, \boldsymbol{\theta}), \ \boldsymbol{D} \tag{3.1a}$$

$$\text{Minimize } C = C(\boldsymbol{u}(\boldsymbol{D}), \boldsymbol{D}) = \int_{-t/2}^{t/2} \int_\omega (\boldsymbol{F} \cdot \boldsymbol{U}) \sqrt{g} \, d\xi^1 d\xi^2 d\xi^3, \tag{3.1b}$$

s.t.

$$\int_{-t/2}^{t/2} \int_\omega \mathrm{H}(\phi_\omega^s(\xi; \boldsymbol{D})) \left(\mathcal{C}^{\alpha\beta\lambda\mu} e_{\alpha\beta}(\boldsymbol{U}) e_{\lambda\mu}(\boldsymbol{V}) + \mathcal{D}^{\alpha\lambda} e_{\alpha 3}(\boldsymbol{U}) e_{\lambda 3}(\boldsymbol{V})\right) \sqrt{g} \, d\xi^1 d\xi^2 d\xi^3$$
$$= \int_{-t/2}^{t/2} \int_\omega (\boldsymbol{F} \cdot \boldsymbol{V}) \sqrt{g} \, d\xi^1 d\xi^2 d\xi^3, \quad \forall \boldsymbol{V} = \boldsymbol{V}(\boldsymbol{v}, \boldsymbol{\eta}) \in U_{ad}, \tag{3.1c}$$

$$\int_{-t/2}^{t/2} \int_\omega \mathrm{H}(\phi_\omega^s(\xi; \boldsymbol{D})) \sqrt{g} \, d\xi^1 d\xi^2 d\xi^3 \leq \overline{V} \int_{-\frac{t}{2}}^{\frac{t}{2}} \int_\omega \sqrt{g} \, d\xi^1 d\xi^2 d\xi^3, \tag{3.1d}$$

$$\boldsymbol{U} \in U, \boldsymbol{D} \in \mathcal{U}_D. \tag{3.1e}$$

In Eq. (3.1), the primary displacement $\boldsymbol{U} = \boldsymbol{U}(\boldsymbol{u}, \boldsymbol{\theta})$ belonging to a prescribed constraint set $U$ and the virtual displacement $\boldsymbol{V} = \boldsymbol{V}(\boldsymbol{v}, \boldsymbol{\eta})$ belonging to an admissible set $U_{ad}$ are assumed to have the following forms [41]:

$$U = U(\xi^1, \xi^2, \xi^3) = u(\xi^1, \xi^2) + \xi^3 \theta_\lambda(\xi^1, \xi^2) a^\lambda(\xi^1, \xi^2), \qquad (3.2a)$$

$$V = V(\xi^1, \xi^2, \xi^3) = v(\xi^1, \xi^2) + \xi^3 \eta_\lambda(\xi^1, \xi^2) a^\lambda(\xi^1, \xi^2), \qquad (3.2b)$$

where $a^\lambda$, $\lambda = 1, 2, 3$ are the corresponding contravariant base vectors associated with the three covariant base vectors $a_\alpha = \partial \varphi / \partial \xi^\alpha$, $\alpha = 1, 2$ and $a_3 = (a_1 \times a_2)/\|a_1 \times a_2\|$, respectively. In Eq. (3.1b) and Eq. (3.1c), $F = F(\xi^1, \xi^2, \xi^3)$ is the prescribed external load, $\mathcal{C} = (\mathcal{C}^{\alpha\beta\lambda\mu})$ and $\mathcal{D} = (\mathcal{D}^{\alpha\lambda})$ are the modified constitutive tensors of the linear elastic material expressed in the convected curvilinear coordinate system $(\xi^1, \xi^2, \xi^3)$. In addition, $e_{\alpha\beta}(U)$ ($e_{\alpha\beta}(V)$) and $e_{\alpha 3}(U)$ ($e_{\alpha 3}(V)$) are the surface and transverse shear strain tensors corresponding to the primary (virtual) displacement, respectively. The specific forms of $e_{\alpha\beta}(U)$ and $e_{\alpha 3}(U)$ depend on the shell theory (constructed based on different stress and/or strain assumptions) adopted. Furthermore, the quantity $\sqrt{g}$ in Eq. (3.1c) and Eq. (3.1d) takes the form of $\sqrt{g} = |g_1 \cdot (g_2 \times g_3)|$ with $g_\alpha = \partial \overline{\varphi}/\partial \xi^\alpha$, $\alpha = 1, 2, 3$. We refer the readers to [41] for more details on the variational formulation for shell analysis. It is also noting that topology optimization of a shell structure can be achieved by finding the material distribution on its mid-surface, which can be characterized by the TDF $\phi_\omega^s(\xi; D)$ in Eq. (3.1) defined on the $\omega$. In addition, $H = H(\cdot)$ is the Heaviside function and $\mathcal{U}_D$ is the set that $D$ belongs to.

**3.2 Structural response and sensitivity analysis**

In the present work, the structural response is calculated approximately using the S3 element (a conventional stress/displacement shell element with three nodes) provided in ABAQUS, which is constructed from a refined shell theory [42]. In order to establish a seamless link with the CAD modeling approaches, triangulated unstructured meshes are adopted for finite element discretization. For the sake of computation efficiency, the ersatz material model is adopted to calculate the element stiffness matrix [43]. For the $e$-th finite element, its equivalent Young's modulus can be calculated in terms of the element equivalent density $\rho_e$ and Young's modulus of the solid material $E^s$ as

$$E_e = \rho_e E^s, \qquad (3.3)$$

where

$$\rho_e = \sum_{j=1}^{3} \mathrm{H}_{\alpha,\epsilon}\left(\left(\phi_{s_\Delta}^s\right)_{e,j}\right)/3 \qquad (3.4)$$

with

$$H_{\alpha,\epsilon}(x) = \begin{cases} 1, & \text{if } x > \epsilon, \\ \dfrac{3(1-\alpha)}{4}\left(\dfrac{x}{\epsilon} - \dfrac{x^3}{3\epsilon^3}\right) + \dfrac{1+\epsilon}{2}, & \text{if } |x| \leq \epsilon, \\ \alpha, & \text{otherwise} \end{cases} \quad (3.5)$$

denoting the regularized Heaviside function ($\epsilon = 0.1$ and $\alpha = 10^{-3}$, respectively, in the present work).

It is worth noting that although the variational formulation for structural analysis in Eq. (3.1) is expressed in the parametric domain for describing a more general problem setting, its numerical implementation is actually achieved by finite element discretization in the physical domain. Since the triangulated mesh for finite element analysis (FEA) is generated on $S_\Delta$, which is an approximation of $S$ (the mid-surface of the shell) in the physical domain, $\phi^s_{S_\Delta} = \phi^s_{S_\Delta}(x)$ (calculated from $\phi^s_{\mathcal{M}} = \phi^s_{\mathcal{M}}(p)$ defined on the parametric domain with the uses of the conformal mapping technique described in Section 2) should be used to characterize the material distribution on $S_\Delta$. In Eq. (3.2), $\left(\phi^s_{S_\Delta}\right)_{e,j}$ denotes the value of $\phi^s_{S_\Delta}$ on the $j$-th node of the $e$-th element on $S_\Delta$.

For a general objective/constraint function $I$, its variation (i.e., $\delta I$) with respect to the variation of a typical design variable $d$ (i.e., $\delta d$) in the following general continuum setting form (provided that some smoothness conditions on $I$ and regularity requirements on the design domain are satisfied):

$$\delta I = \int_{-\frac{t}{2}}^{\frac{t}{2}} \int_S r(U(x), W(x)) \delta\phi^s_S(x; \delta d) \, dx^1 dx^2 dx^3, \quad (3.6)$$

where $r = r(U(x), W(x))$ is a function of $U(x)$ and $W(x)$ which are the primary and adjoint displacement fields described in the physical domain (for the considered compliance minimization problem $W(x) = -U(x)$), respectively, while $\delta\phi^s_S(x; \delta d)$ denotes the variation of $\phi^s_S = \phi^s_S(x; D)$ due to the variation of $d$ (i.e., $\delta d$). Since $\phi^s_S = \phi^s_S(x; D) = \phi^s_{\mathcal{M}}(p; D)$ with $p = f(x)$ (noting that $f$ is the conformal mapping from $S$ to $\mathcal{M}$). Under this circumstance, we have $\delta\phi^s_S(x; \delta d) = \delta\phi^s_S(p; \delta d) = \delta\phi^s_S(f(x); \delta d)$. Considering the fact that $\mathcal{M} = \bigcup_{k=1}^{N_u} \mathcal{M}_k$, it yields that

$$\delta\phi_{\mathcal{S}}^{\mathrm{S}}(f(x);\delta d)$$

$$= \begin{cases} \left(\dfrac{\partial \phi_{\mathcal{M}_k}^{\mathrm{S}}(\boldsymbol{p};\boldsymbol{D})}{\partial d}\right)\delta d, & \text{if } \boldsymbol{p} = f(x) \in \mathcal{M}_k \setminus \bigcup_{l=1,l\neq k}^{N_{\mathcal{U}}} \mathcal{M}_l, \quad (3.7\mathrm{a}) \\ \left(\dfrac{\partial \mathcal{KS}\left(\phi_{\mathcal{M}_k}^{\mathrm{S}}(\boldsymbol{p};\boldsymbol{D}),\dots,\phi_{\mathcal{M}_m}^{\mathrm{S}}(\boldsymbol{p};\boldsymbol{D}),\right)}{\partial d}\right)\delta d, & \text{if } \boldsymbol{p} = f(x) \in \mathcal{M}_k \cap \dots \cap \mathcal{M}_m. \quad (3.7\mathrm{b}) \end{cases}$$

The calculation of $\partial\phi_{\mathcal{M}_k}^{\mathrm{S}}(\boldsymbol{p};\boldsymbol{D})/\partial d$ follows exactly the same way as that in the MMC approach developed for flat 2D case.

When finite element method is used for approximated structural analysis, the sensitivity of the concerned objective function can be calculated in the following discrete form:

$$\frac{\partial I}{\partial d} = -\boldsymbol{U}^{\mathrm{T}}\frac{\partial \boldsymbol{K}}{\partial d}\boldsymbol{U}, \qquad (3.8)$$

where $\boldsymbol{U}$ is the vector of the displacement field on $\mathcal{S}_\Delta$ and $\partial \boldsymbol{K}/\partial d$ can be determined by the variation of $\delta\phi_{\mathcal{S}_\Delta}^{\mathrm{S}}$ with respect to $\delta d$ in the way described above. Furthermore, the sensitivity of the shell volume with respect to $d$ is quite straightforward and will not be discussed here.

### 3.3 The flowchart of the solution procedure

As a summary, the flowchart of the proposed integrated solution procedure for topology optimization on complex surfaces with arbitrary genus is described in this subsection.

**Step 1: Surface pre-processing**

(a) Generating a triangular surface on the mid-surface $\mathcal{S}$ of a given object $\mathcal{O}$ (whose geometry information can be obtained from CAD modeling or direct 3D scanning) and obtain a surface $\mathcal{S}_\Delta$ constituted by the generated triangulated surfaces as an approximation of $\mathcal{S}$ with enough accuracy.

(b) Dividing $\mathcal{S}_\Delta$ into several parts $\mathcal{U}_k$, $k = 1, \dots, N_\mathcal{U}$ based on its geometric features such that $\mathcal{S} = \bigcup_{k=1}^{N_\mathcal{U}} \mathcal{U}_k$.

(c) Defining the cutting lines $\mathcal{C}_l$, $l = 1, \dots, N_\mathcal{U}$ which make every $\mathcal{U}_k$ ($k = 1, \dots, N_\mathcal{U}$) expand into a single simply connected open surface $\mathcal{U}_k^*$ with genus zero.

**Step 2: Parameterization based on conformal mapping**

a) Looping from $k = 1$ to $N_\mathcal{U}$ for every $\mathcal{U}_k^*$, parameterizing each $\mathcal{U}_k^*$ through the CCM

technique described in Section 2 and obtaining the corresponding conformal mappings $f_k: \mathcal{U}_k^* \to \mathcal{M}_k, k = 1$ to $N_\mathcal{U}$.

**Step 3: Topology optimization**

a) Placing components in each parametric domain $\mathcal{M}_k$, $k = 1$ to $N_\mathcal{U}$.

b) Computing $\phi_{\mathcal{M}_k}^S(\boldsymbol{p})$ based on the design variables associated with each component on every parametric domain $\mathcal{M}_k$ and obtaining $\phi_{\mathcal{S}_k}^S(\boldsymbol{x})$ through $\phi_{\mathcal{S}_k}^S(\boldsymbol{x}) = \phi_{\mathcal{M}_k}^S(f_k(\boldsymbol{x}))$.

c) Obtaining $\phi_{\mathcal{S}}^S(\boldsymbol{x})$ from $\phi_{\mathcal{S}_k}^S(\boldsymbol{x})$, $k = 1$ to $N_\mathcal{U}$ through Eqs. (2.15)-(2.17) in Section 2.

d) Performing topology optimization on $\mathcal{S}_\Delta$ using $\phi_{\mathcal{S}}^S(\boldsymbol{x})$ based on the traditional MMC approach.

## 4. Numerical Examples

In this section, three numerical examples are examined to demonstrate the effectiveness of the proposed framework for topology optimization on surface with arbitrary genus. Triangular meshes are adopted to represent the geometry of complex surfaces. The same mesh is also used for finite element analysis. Unit thickness three-node bilinear shell finite elements are used to solve the structural response. Without loss of generality, all involved quantities are assumed to be dimensionless and the thickness $t$ of the considered shell is $t = 1$. The Young's modulus and the Poisson's ratio of the isotropic solid material are chosen as $E^s = 1$ and $v = 0.3$, respectively. In all examples, the available volume of the solid material is $\bar{V} = 0.4 V_\mathrm{D}$ with $V_\mathrm{D}$ denoting the volume of the design domain on the surface. The Method of Moving Asymptotes [44] is used as numerical optimizer. The optimization process is terminated if the relative change of each design variable Tol between two consecutive iterations is below a specified threshold (i.e., Tol=0.001).

### 4.1 Saddle-shaped shell example

In this example, we consider a shell with saddle-shaped mid-surface with its geometry and boundary conditions shown in Fig. 7. The saddle-shaped shell structure is subjected to a horizontal tangential concentrated load at its saddle point. Although this problem is symmetric, the entire structure is optimized to test the robustness of the proposed approach. The design domain on the saddle-shaped mid-surface of the shell is discretized into 29216 triangular meshes with 14868 nodes for geometric description and finite element analysis.

From topology point of view, the genus of an open saddle surface is zero and is globally homeomorphic to a planar rectangle. This means that the corresponding computational conformal mapping function can be determined without any cutting operation (see Fig. 8). As shown in Fig. 9, there are totally 16 components (containing $7 \times 16 = 112$ design variables) are distributed in the planar rectangle parametric domain characterizing the topology of the initial design. Fig. 9 shows the corresponding component distribution in the physical domain (the mid-surface of the shell). Fig. 10 plots the iteration history of optimization process. It is found that the structural compliance experiences a rapid drop in the first 5 iterations and then begins to decrease gradually in the following steps and finally converges to $I^{\mathrm{opt}} = 22.56$ at about the 120th step. Some intermediate optimization results are also presented in Fig. 10. It can be observed that the distributions of components in the parameter domain and the physical domain do maintain the topological consistency. As the optimization iterations proceed, the components gradually form a connected load transmission path between the loading point and the fixed structural boundary.

The final optimization results are plotted in Fig. 11 and it can be seen that the corresponding components form a clear load transmission path on the mid-surface of the shell. Although no symmetry constraints are imposed on the optimization problem, the optimized structure still maintains the symmetry exactly. It is also found that the stress is higher in the region near the saddle point in the initial design due to the existence of concentrated horizontal forces. While in the final design, the optimized components form an elliptical region around the saddle point automatically, which effectively relieves the stress concentration phenomenon.

**4.2 Torus-shaped shell example**

In this example, topological design of a torus-shaped shell is considered to demonstrate the capability of the present approach to deal with surface of non-zero genus and its potential of being integrate seamlessly with 3D scanning technique. The problem under consideration is shown in Fig. 12. The inner ring of the shell assumed to be fixed and four rotationally symmetric shear forces are applied at four points along the outer ring.

Actually, with the help of modern 3D scanning technique, a high-precision discrete point cloud obtained from the scanning of a surface can be generated efficiently and the corresponding data representing the geometry of the surface can be exported in a standard PLY (Polygon) format [45]. Without resorting to any further post-processing steps, the vertex and face information contained in the exported PLY data can be used directly to generate the required triangular mesh for geometric description and finite element analysis in the proposed method. Fig. 13 shows the discrete geometry model of the considered torus surface structure generated by 3D scanning with 31840 vertices and 63680 triangular facets. Since the genus of the torus is non-zero, it cannot be mapped conformally to a single rectangular in the parametric domain. Under this circumstance, as described in Section 2 and shown in Fig. 13, we first cut the torus along the path indicated by two intersecting circles and then map the intermediate surface obtained by cutting operation to a rectangular. It is worth mentioning that, actually, as the genus of a surface increases, the cutting path is hard to be determined by intuition. Fortunately, general algorithms have already been developed to determine the cutting path automatically for surfaces with arbitrary genus in the field of computational topology [46, 47] and can be used to construct conformal mapping for complex surfaces.

Topology optimization can be performed once the mapping between the intermediate surface obtained by cutting operation and a planar rectangular has been established. Fig. 14 shows the initial design containing 64 components distributed in the rectangle and their image on the torus surface, respectively. The iteration history of the optimization process and the optimized designs are provided in Fig. 15 and Fig. 16, respectively. An optimized design with $I^{\mathrm{opt}} = 1.80$ is obtained after 320 steps. It can be observed from these figures that driven by the optimization algorithm, the components automatically achieve a smooth connection along the cutting boundary even though no special constraints are introduced. A stable structural topology is achieved after 80 optimization iteration steps while the subsequent steps only adjust some minor structural details. In the optimized design, the components form a lattice-like structure on the torus surface, which is believed to be very efficient in resisting torsional deformation [48].

### 4.3 Tee-branch-shaped shell example

In this final example, we consider a complex tee-branch pipe structure which can be modeled as a thin shell to illustrate the effectiveness of the proposed multi-patch stitching technique. The geometry of the tee-branch pipe, external load, and boundary conditions are all shown in Fig. 17. Noting that not only the structure has complex topology, but each branch of the pipe also has different shape of non-uniform cross-section as indicated in the Figure. Although the cutting operation can also be applied to establish a global conformal mapping, the topology complexity of the tee-branch pipe may inevitably cause excessive distortions of the structural components on the surface although the components distributed in the parametric domain are of regular shapes. Therefore, based on the intrinsic topology character of the tee-branch structure, we first partition it into four patches (including three branch patches and a joint patch) as shown in Fig. 18 (a) and then computational conformal mapping technique is used to establish the surface parameterization of each patch. Each branch patch is actually topologically equivalent to a cylindrical surface and can be directly mapped to a rectangular by cutting along the direction of the generator. Therefore, the corresponding conformal mappings can be established in an ordinary way. However, after cutting along the selected lines, it is still impossible to establish a topological homeomorphism between the joint patch and a rectangular plane due to the existence of a "hole" on this patch, as shown in (Fig. 18(c)). In other words, the joint patch is topologically equivalent to a rectangular with a hole. To tackle this problem, one approach is to make another cut and turn the joint patch into a simply-connected open surface. Here, we, however, solve this in another way by first filling the "hole" in the joint patch by Delaunay triangularization [49, 50], making it topologically equivalent to a rectangular, establishing the conformal mapping and finally deleting the triangular mesh corresponding to the "hole" in the rectangular (Fig. 18(c)). With the use of the above treatment, a conformal mapping relationship between the joint patch and a rectangular with a hole can be established. Once the corresponding conformal mappings are constructed for the four patches, the global TDF characterizing the material distribution on the tee-branch pipe can be determined in the way described in Section 2.

Fig. 19 shows the initial component layout and the process of assembling the components from the different patches. Fig. 20 shows the iteration history of the optimization process and it can be observed that an optimized design can be found after 200 iteration steps. In the optimized design ($I^{\mathrm{opt}} = 141.04$) shown in Fig. 21, the material distributions in the two lower branches of the pipe are very similar to that in the well-known 2D short beam example [16], while the overall structural topology of the optimized tee-branch pipe has some similarity to that of the optimized structural layout of the planar double L-bracket structure [51]. Fig. 21 plots the final component distribution on each patch and illustrates the assembly process of the optimization results on each patch.

## 5. Concluding remarks

In the present work, an integrated paradigm for topology optimization on complex surfaces with arbitrary genus is proposed. The two supporting pillars are the Moving Morphable Component (MMC)-based approach, where the topology of a structure can be described by a set of parameters explicitly, and the computational conformal mapping technique with which a complex surface with arbitrary genus can be mapped into a set of regular 2D parameter domains in a systematic way numerically. The effectiveness and the applicability of the proposed paradigm for topological design on complex 2D surfaces have also been verified by several numerical examples provided.

The advantages of the proposed paradigm can be summarized from the following aspects: (1) The proposed solution paradigm is actually based on unstructured triangularization of surface. This treatment naturally renders its applications to a variety of complex surfaces constructed from different ways (e.g., analytical description, CAD modeling and 3D-scaning generated point cloud). It can also be used to perform topology optimization directly based on the triangular meshes generated by CAE software. It is also worth noting that since the conformal mapping is only needed to be established once before optimization, it will not introduce too much computational cost. (2) Since MMC method is adopted for topology optimization, the number of design variables can be reduced significantly and clear load transmission paths can be identified easily in final optimized designs. Moreover, since the constructed conformal mapping is topology preserved, the load path identification approach developed in traditional MMC framework [52] can also be used to eliminated inactive degrees of freedom from the finite element model to speed up the finite element analysis. For the limitation of space, this feature will be demonstrated in detail in a separate work. (3) The developed multi-patch stitching technique can greatly enhance the fidelity of geometry modeling and therefore effectively reduce the nonlinearity of the constructed conformal mapping through local mapping assembling compared to the case where a stiff global conformal mapping is established based on a single patch. It is also very helpful for alleviating the mismatch of components deployed on the boundaries of neighboring patches. Furthermore, it is worth noting that the conformal mapping is needed to be established only once before optimization and therefore its computation will not deteriorate the efficiency of the optimization process.

The present work can be extended along various directions. For example, although only topology optimization of mechanical systems is considered in this work, the proposed solution paradigm can be extended to solve more general surface design problems considering multi-physics effects such as heat transfer control and electromagnetic wave guidance. It can also be generalized to tackle the problem of multi-scale design on complex surfaces by combining the MMC-based techniques developed for problems in 2D flat space [53]. Moreover, the proposed paradigm is also applicable to optimize the layout of stiffeners on shells with complex spatial geometries. This can be achieved by modeling the stiffeners as a set of morphable components with explicit geometry descriptions moving on the shells [54]. Corresponding research results will be reported in separate works.

**Acknowledgment**

This work is supported by the National Natural Science Foundation (11821202, 11732004, 12002077, 12002073), the National Key Research and Development Plan (2020YFB1709401), and 111 Project (B14013). Special thanks also go to Prof. Gary Pui-Tung Choi from Massachusetts Institute of Technology for his generous help on understanding the theory of computational conformal mapping and sharing us with the codes for its numerical implementation.

**Figures**

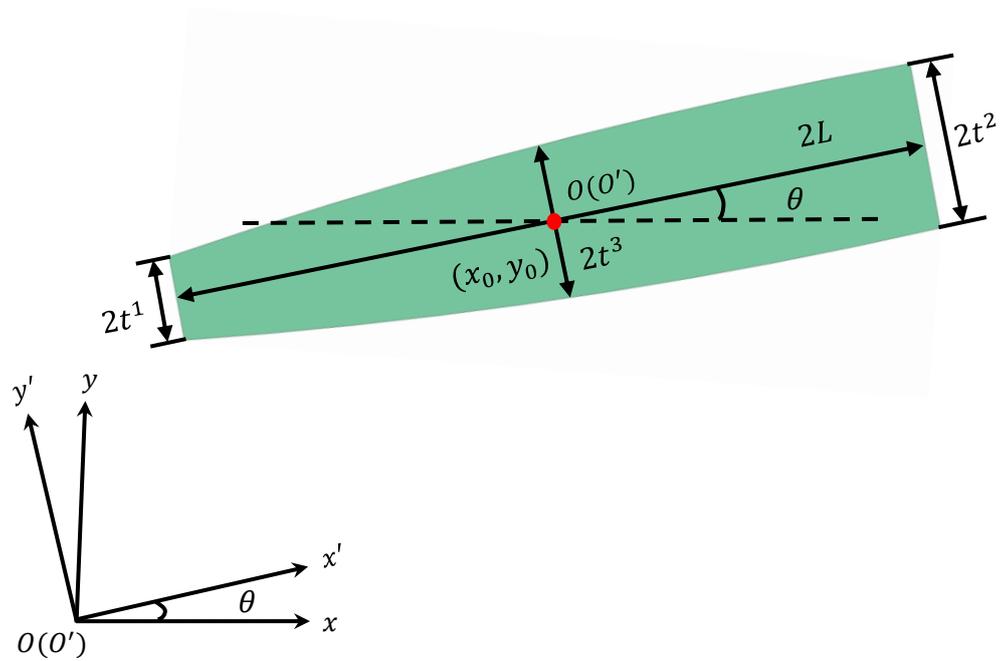

Fig. 1　Geometry description of a typical 2D structural component.

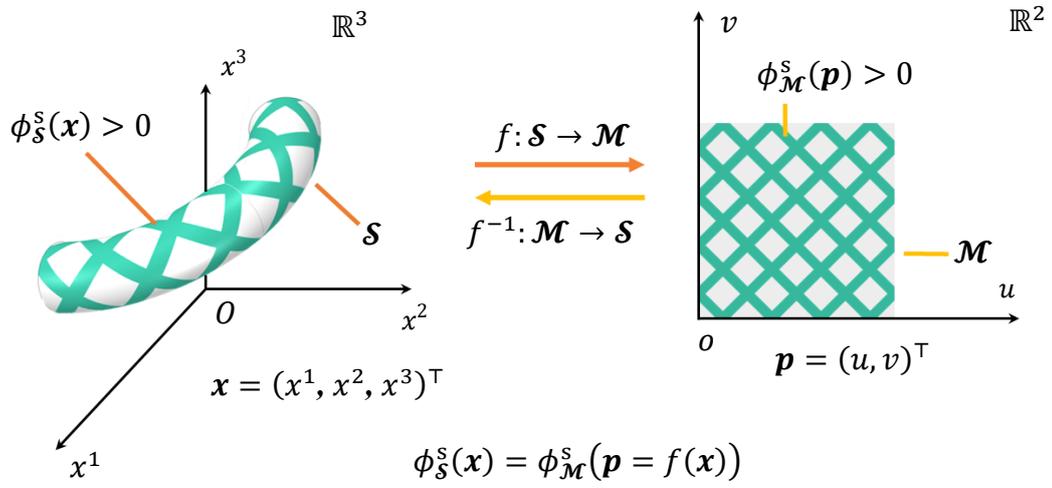

Fig. 2　Parameterization of a surface embedded in 3D Euclidean space and the construction of the corresponding TDF.

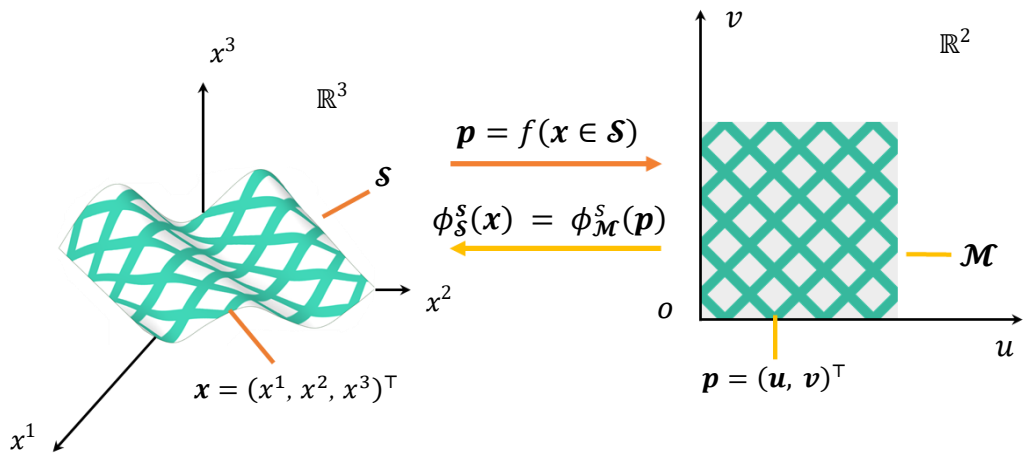

Fig. 3　TDF definition on a simply-connected open surface with genus zero.

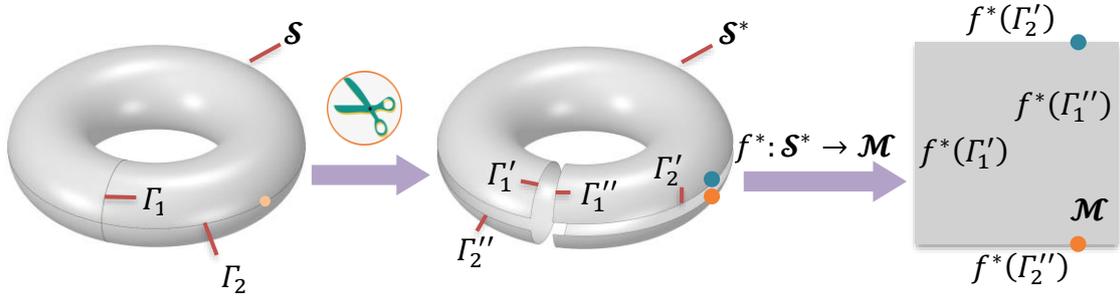

(a) Parameterization of a surface with non-zero genus via cutting operation.

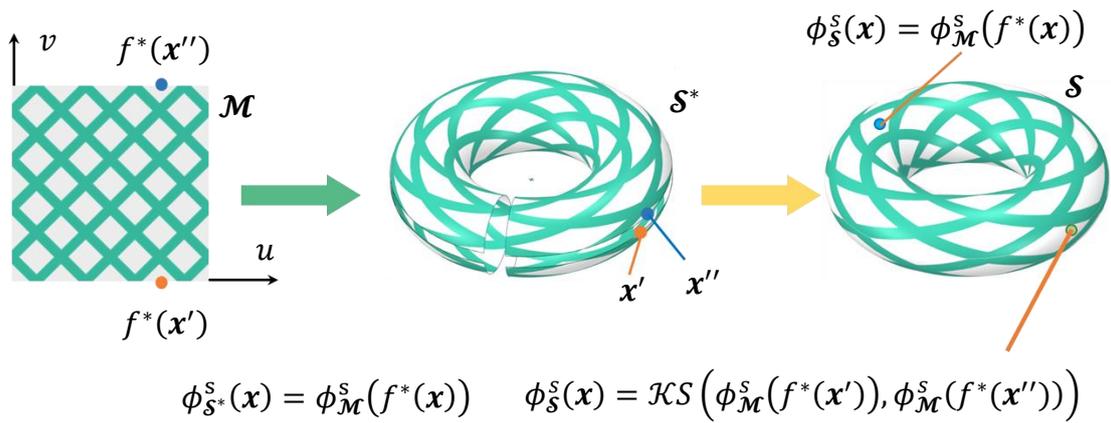

(b) TDF definition on a surface with non-zero genus.

Fig. 4  Parameterization and TDF construction of a surface with non-zero genus.

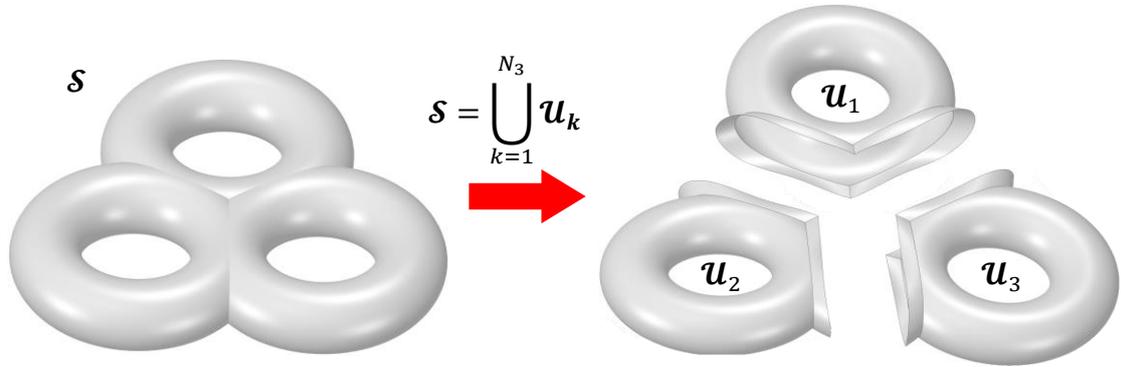

(a) Decomposing a surface into several patches.

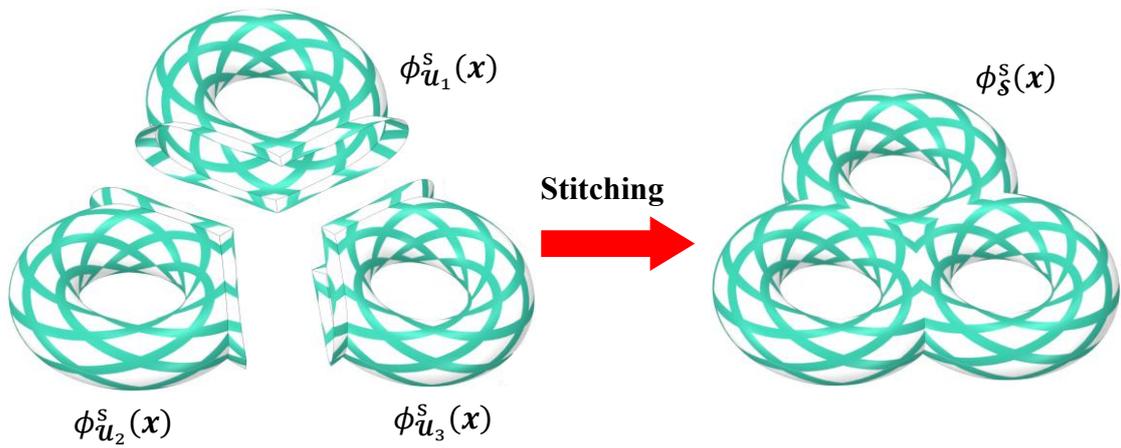

(b) TDF construction on each patch and the corresponding stitching operation.

Fig. 5  TDF construction on a complex surface by multi-patch stitching technique .

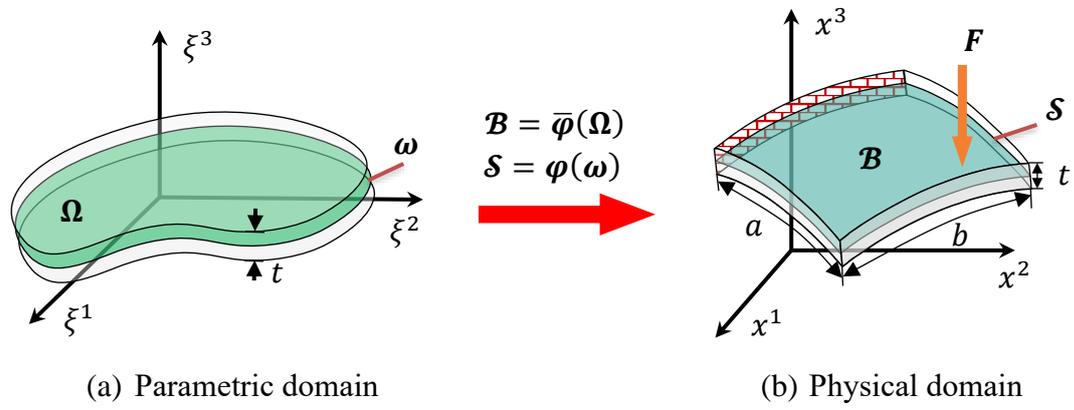

(a) Parametric domain  (b) Physical domain

Fig. 6  Minimizing the compliance of a shell structure by topology optimization.

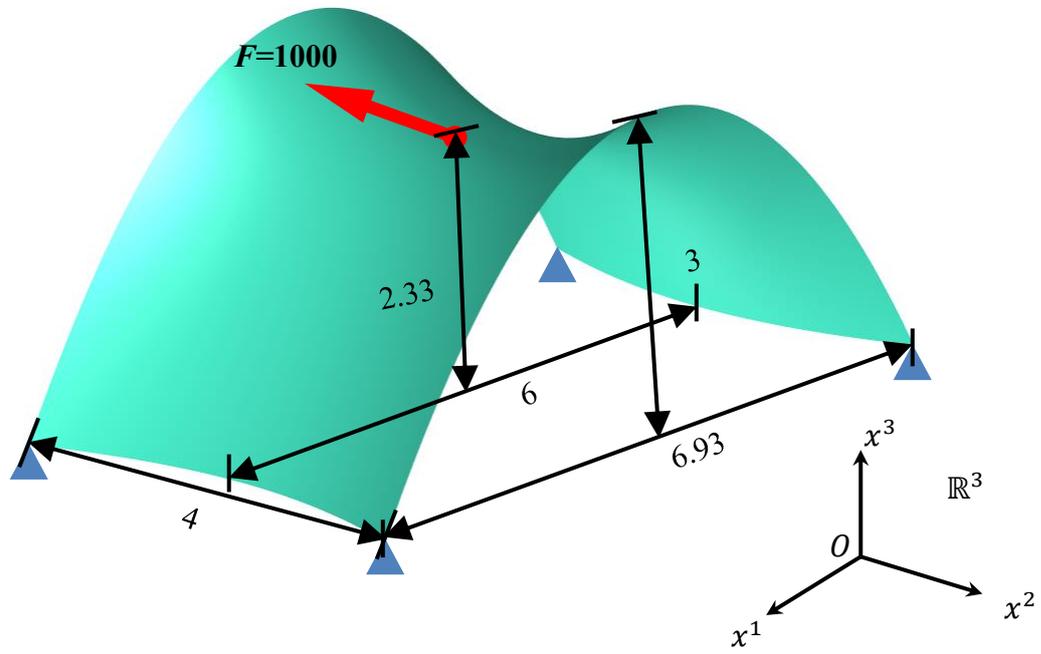

Fig. 7 The problem setting of the saddle-shaped shell example.

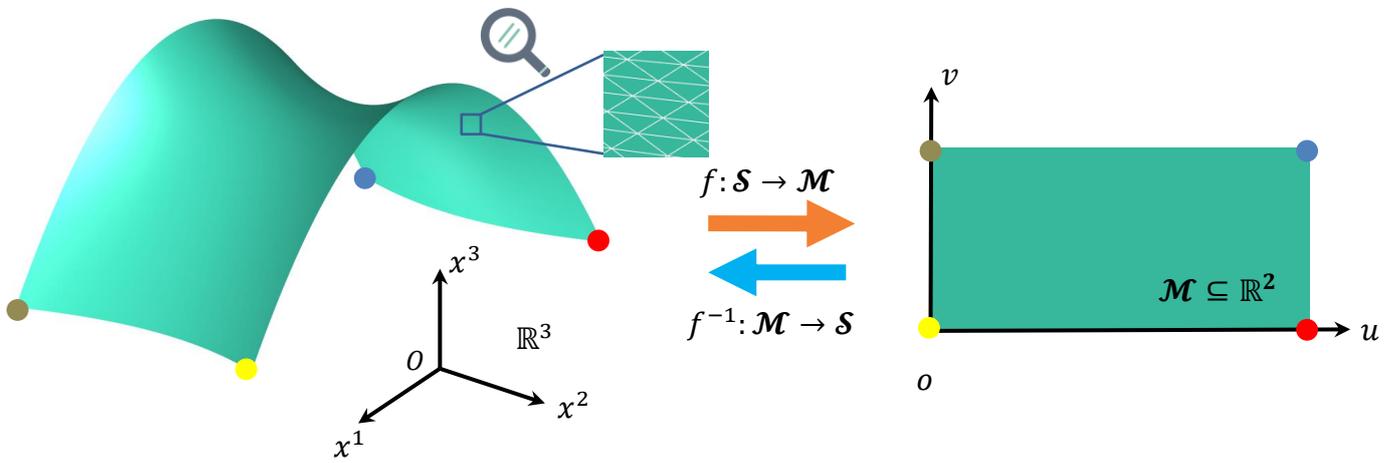

Fig. 8　Parameterization of the mid-surface of the saddle-shaped shell by computational conformal mapping.

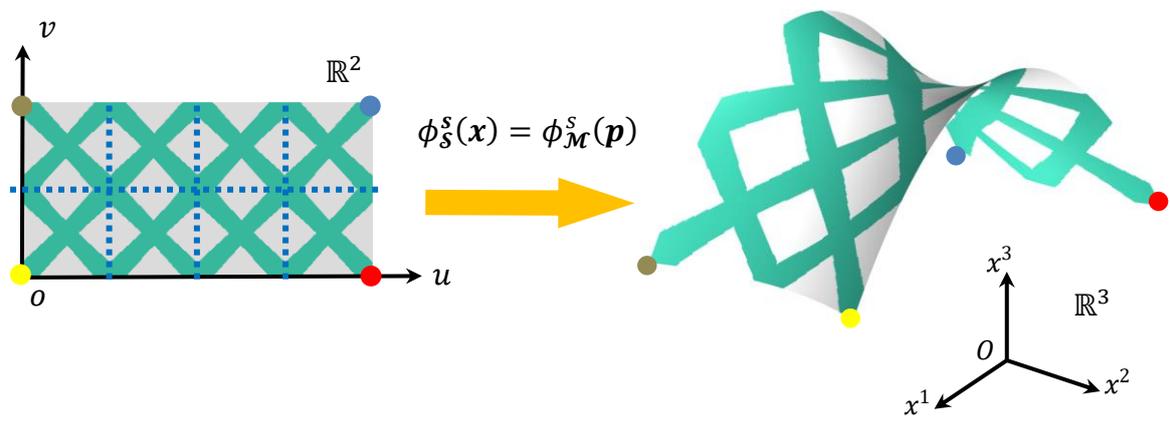

Fig. 9　Initial components layout of the saddle-shaped shell example.

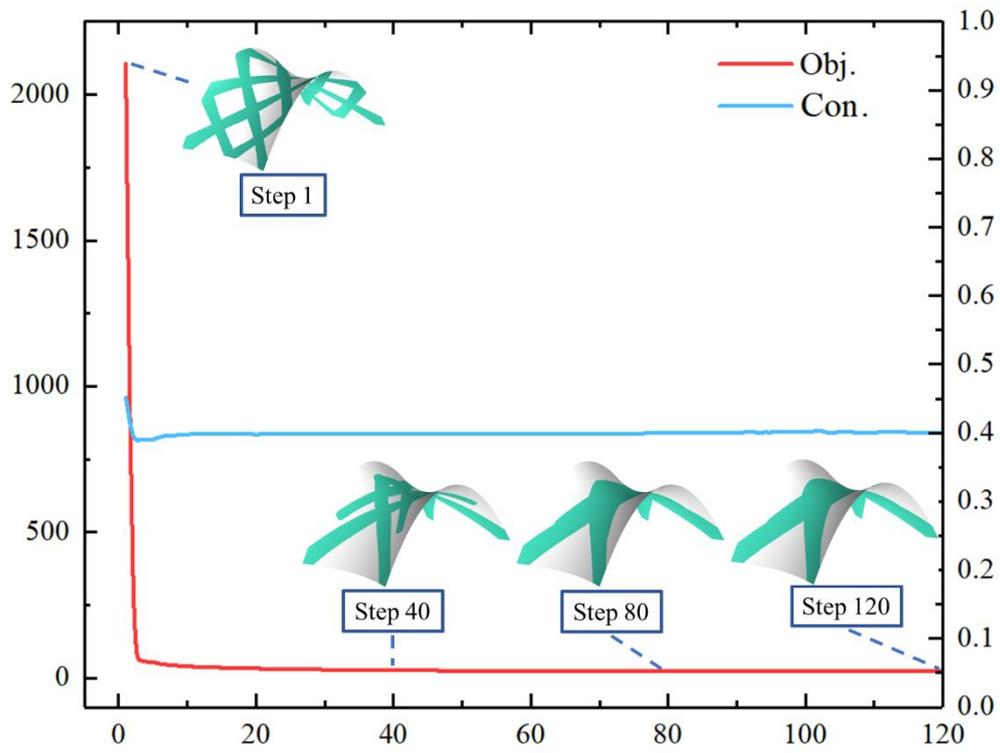

Fig. 10  Iteration history of the saddle-shaped shell example.

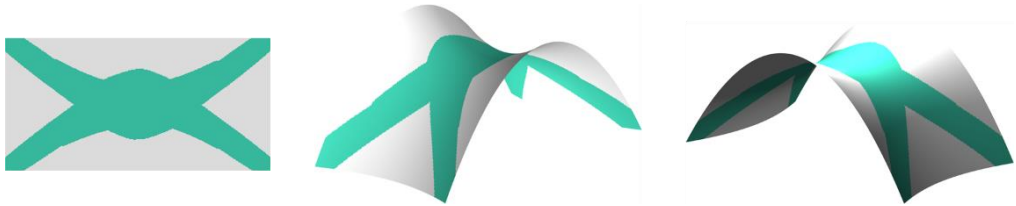

Fig. 11 The optimized structure of the saddle-shaped shell example.

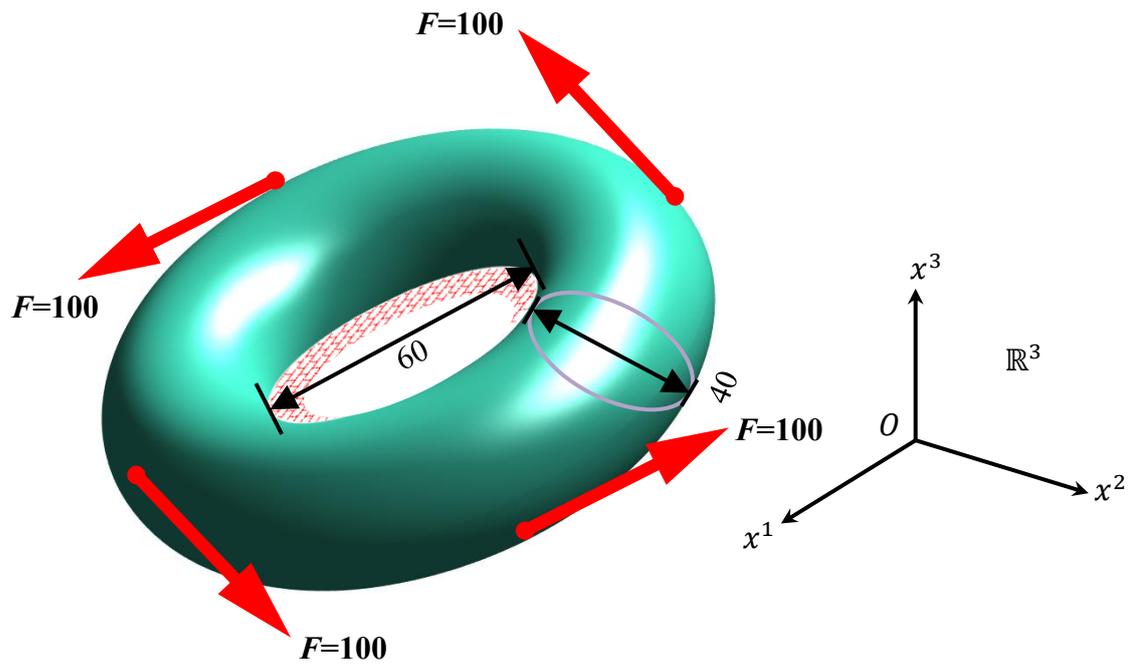

Fig. 12 The problem setting of the torus shell surface.

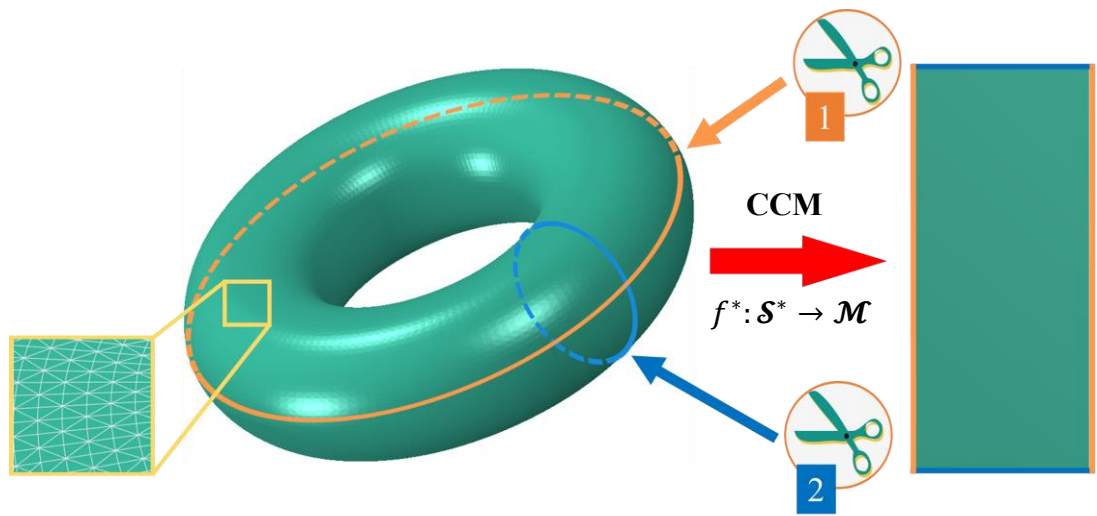

Fig. 13    Torus surface parameterization via cutting operation.

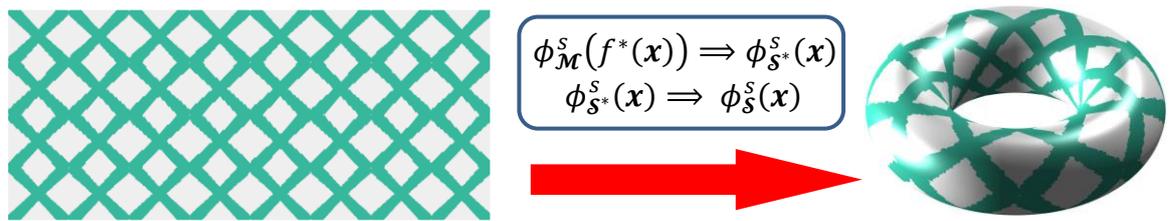

Fig. 14  Initial components layout of the torus-shaped example.

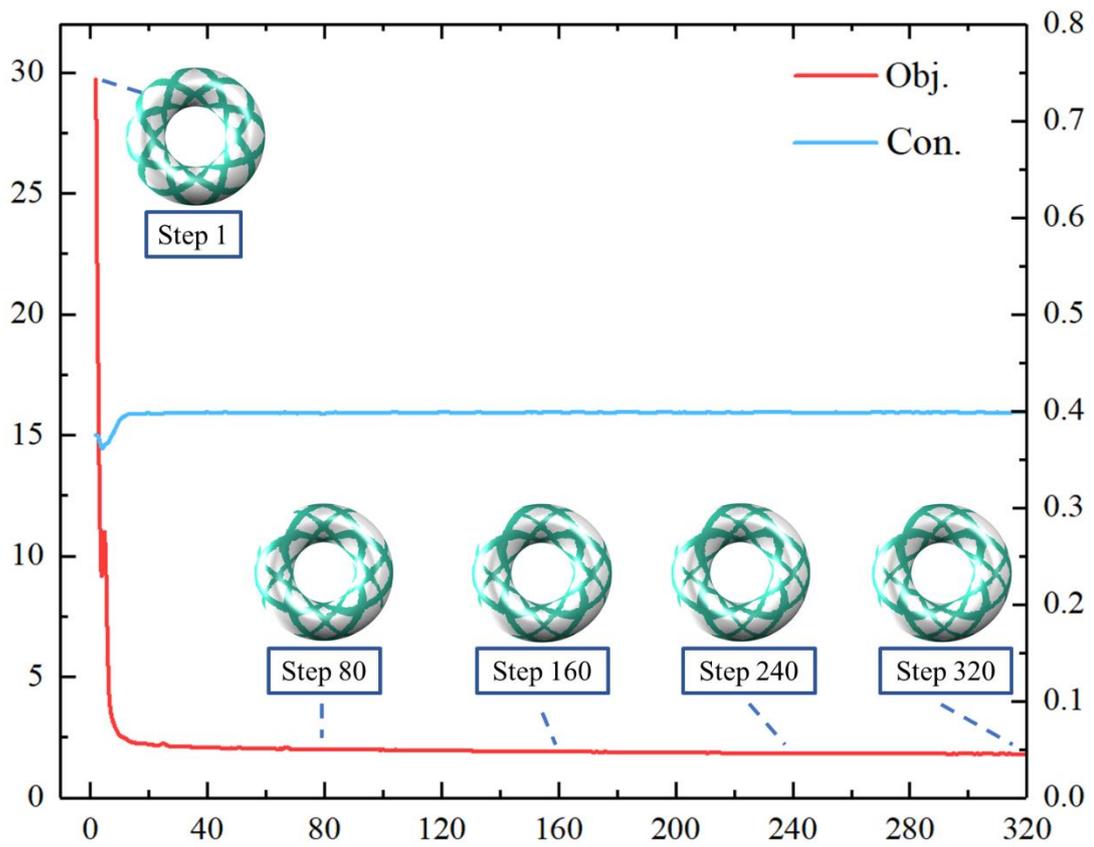

Fig. 15    The iteration history of the torus surface example.

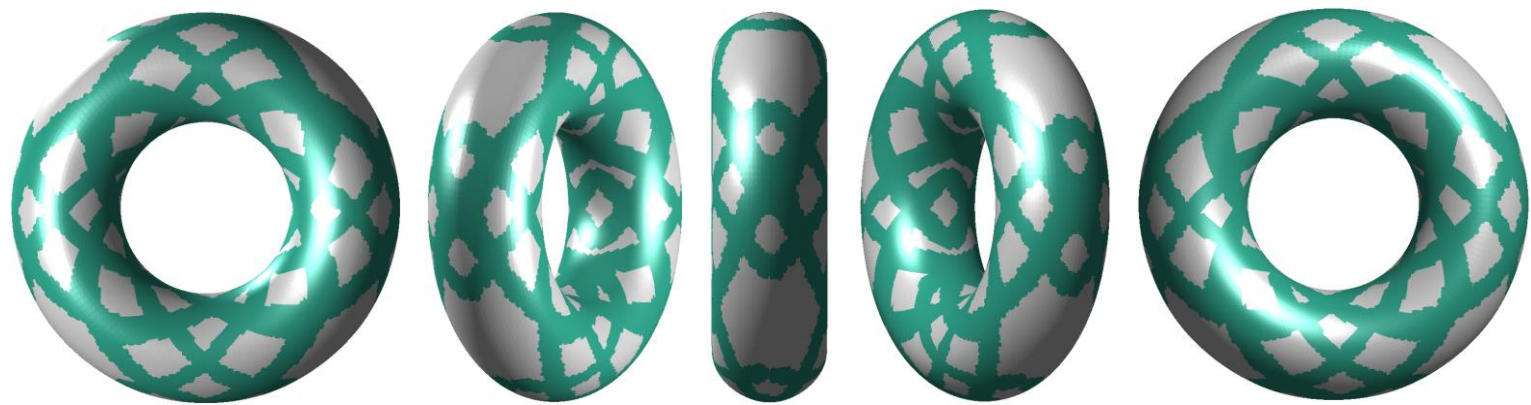

Fig. 16 The optimized design of the tour surface example (viewing from different directions).

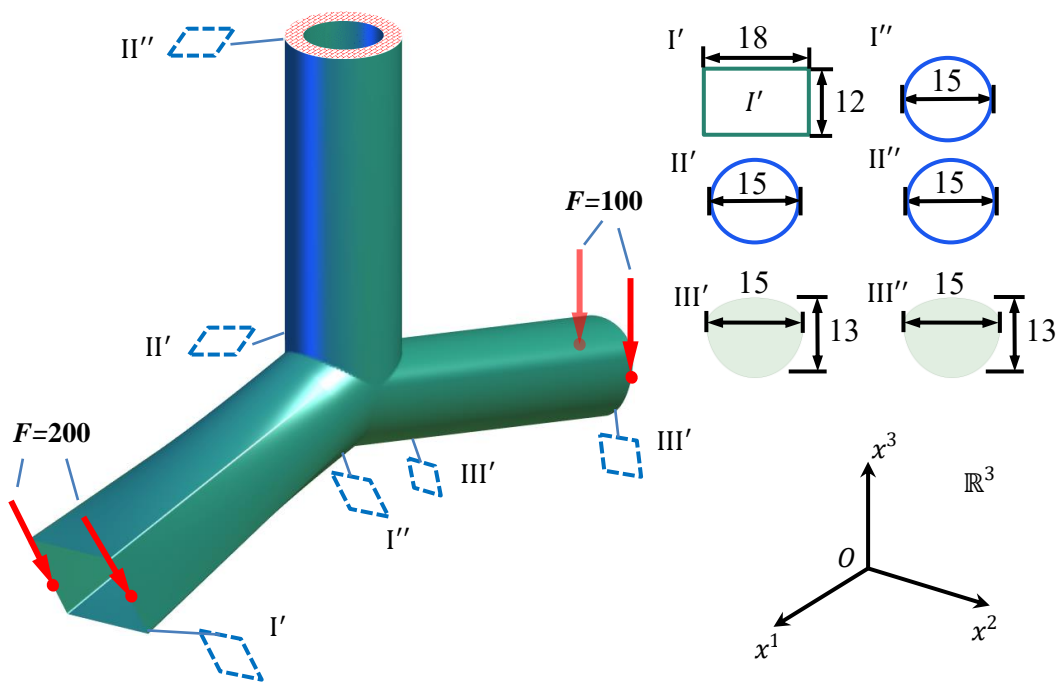

Fig. 17 Problem setting of the tee branch pipe example.

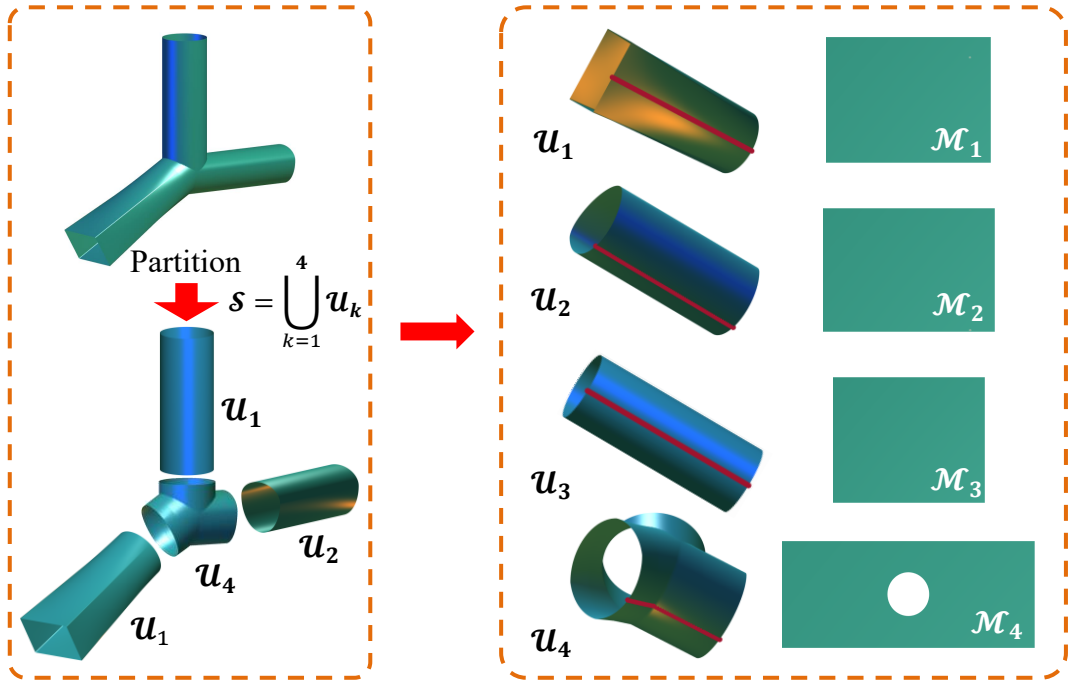

(a) Partition the surface into four patches.

(b) Surface parameterization of each patch.

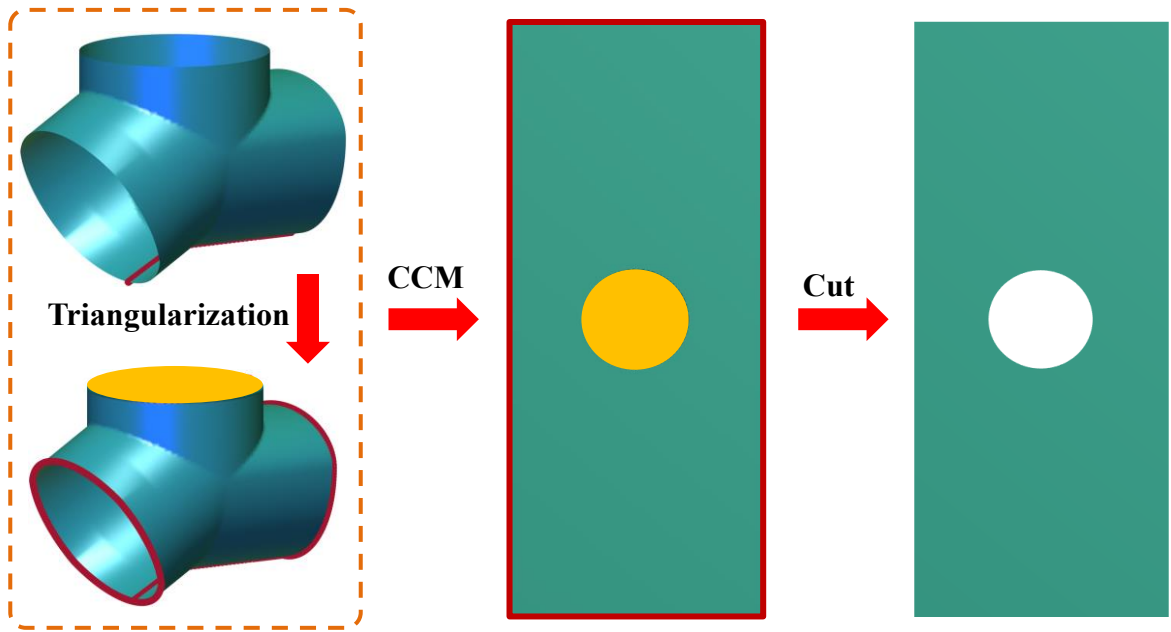

(c) Delaunay triangularization and parameterization.

Fig. 18  Surface parameterization via multi-patch stitching approach.

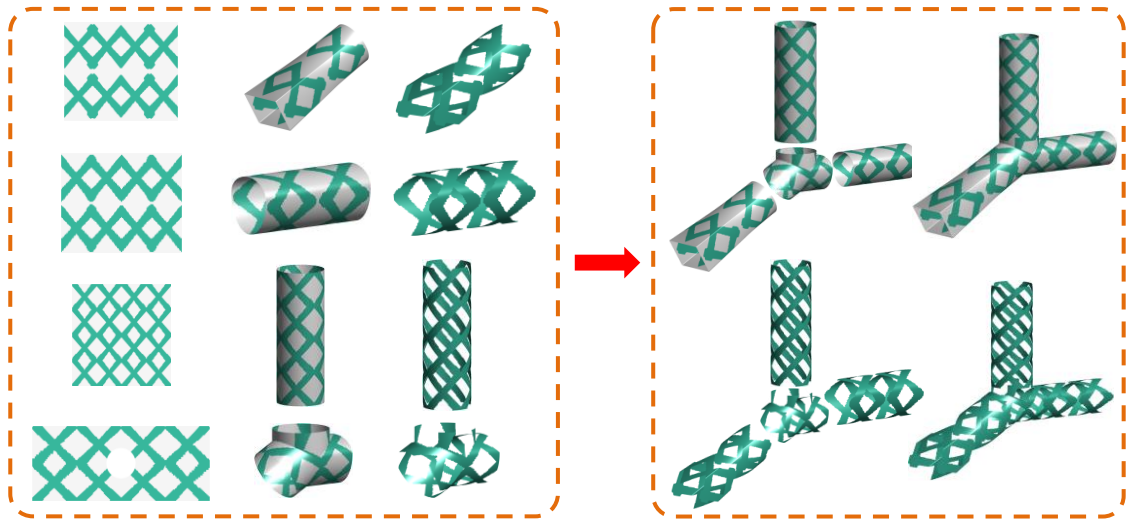

(a) Components on each patches.  (b) Assembling of different patches.

Fig.19　Initial components layout of the tee branch pipe example and the mappings of different patches.

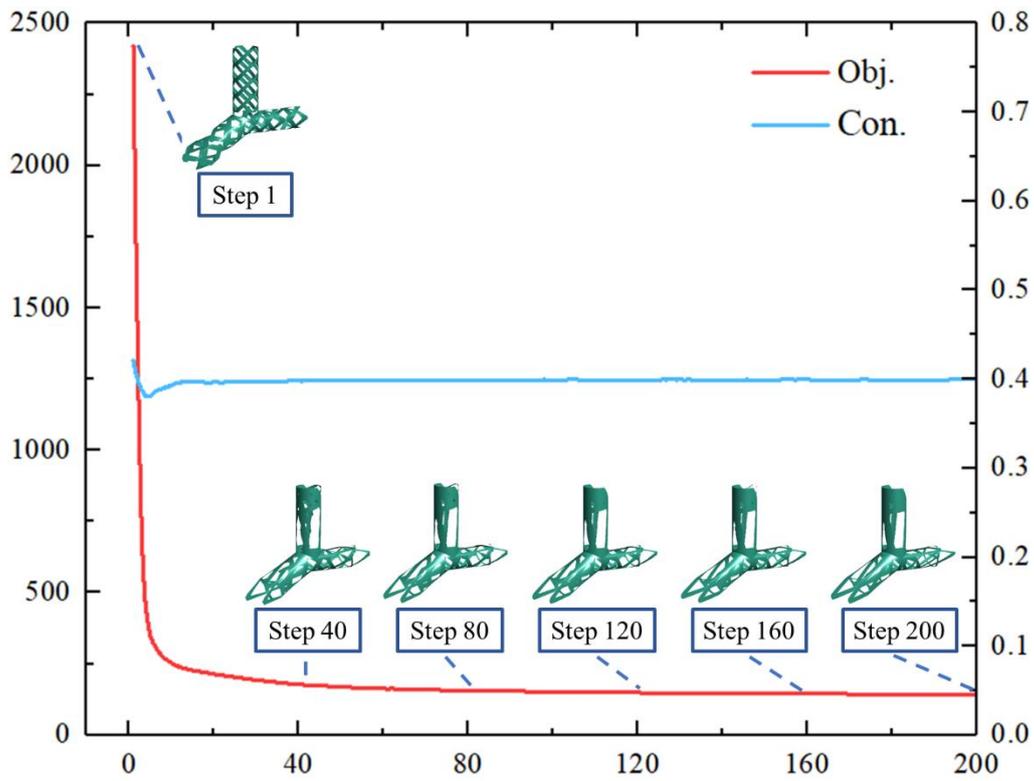

Fig.20　The Iteration history of the tee branch pipe example.

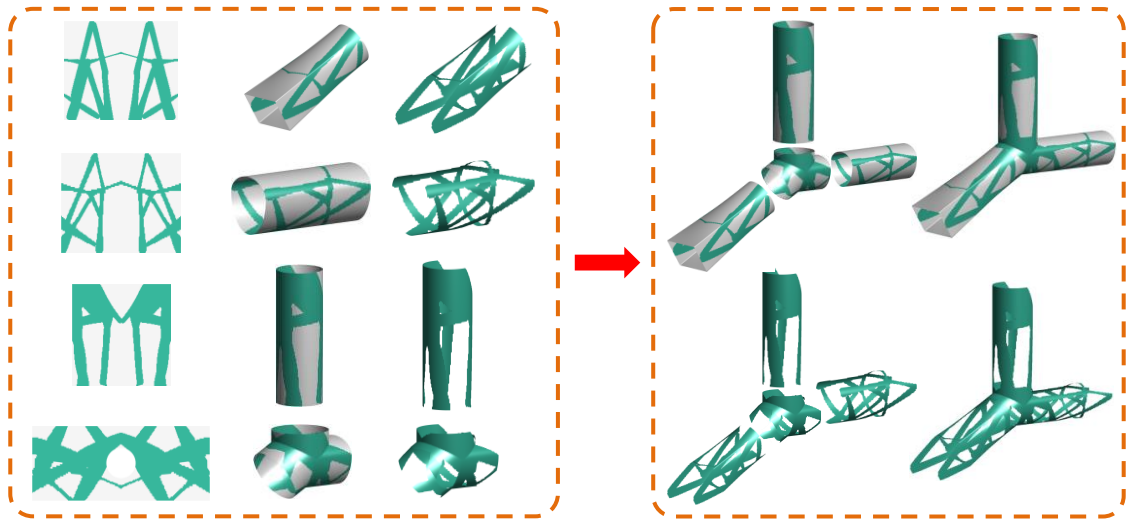

(a) Mapping of individual patch.  (b) Assembly of different patches.

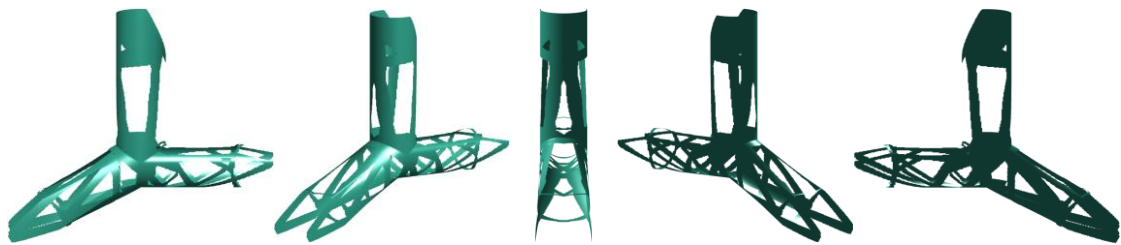

(c) The optimized design (viewing from different directions).

Fig.21  The optimized design of the tee branch pipe example.